# Upper Entropy Axioms and Lower Entropy Axioms

Jin-Li Guo[1], Qi Suo
Business School, University of Shanghai for Science and Technology,
Shanghai 200093, China

The paper suggests the concepts of an upper entropy and a lower entropy. We propose a new axiomatic definition, namely, upper entropy axioms, inspired by axioms of metric spaces, and also formulate lower entropy axioms. We also develop weak upper entropy axioms and weak lower entropy axioms. Their conditions are weaker than those of Shannon-Khinchin axioms and Tsallis axioms, while these conditions are stronger than those of the axiomatics based on the first three Shannon-Khinchin axioms for superstatiscs. The expansibility, subadditivity and strong subadditivity of entropy are obtained in the new axiomatics. Tsallis statistics is a special case of the superstatistics which satisfies our axioms. Moreover, different forms of information measures, such as Shannon entropy, Daroczy entropy, Tsallis entropy and other entropies, can be unified under the same axiomatic framework.

The Boltzmann–Gibbs (BG) statistics works perfectly for classical systems with short range forces and relatively simple dynamics in equilibrium. However, Einstein never accepted Boltzmann's principle $H = k \log N$, because he argued that the statistics ($N$) of a system should follow from its dynamics and, in principle, could not be postulated a priori (1). The foundation of Boltzmann statistical mechanics is the Boltzmann-Gibbs entropy. Entropy plays an important role in statistical mechanics.

Weakly interacting statistical systems are perfectly described by thermodynamics. However, one may be interested in nonequilibrium systems with a stationary state, long transient behaviour in systems with long-range interactions. In a system with long-range interactions, its sub-systems interact with each other. The energy of a system may not simply be the sum of the subsystem energy. Obviously, the entropy of such a system may no longer be the sum of its component entropy. It has been realized that systems with long-range interactions cannot be described by extensive Boltzmann statistics. The statistics sometimes may even lead to some unreasonable results. E.g., from the view of Boltzmann statistics, the long-range interactions may lead to a thermodynamic singularity. In physics there is no way to obtain sum entropy in a finite system. Therefore, the generalized entropy function is regarded as a starting point for universal statistical mechanics.

It was Tsallis in his seminal paper (2) who for the first time suggested to generalize statistical mechanics using his entropic forms. Tsallis, Beck, Cohen et. al. generalized Boltzmann-Gibbs function and established new statistical mechanics, namely superstatistics (1, 2, 3). The conception "statistics of statistics" was firstly put forward by Beck et. al. (1), which seemed to be more

---

effective to describe nonequilibrium systems. Superstatistics describes statistical systems that behave like superpositions of different inverse temperatures $\beta$, so that the probability distribution is $p(\varepsilon_i) \propto \int_0^\infty f(\beta)e^{-\beta\varepsilon_i}d\beta$, where the "kernel" $f(\beta)$ is nonnegative and normalized [$\int_0^\infty f(\beta)d\beta = 1$]. Hanel and Thurner et al.(4) discussed the relation between this distribution and the generalized entropic form $H = \sum_i g(p_i)$. The maximum entropy principle is a method for obtaining the most likely distribution functions of observables from statistical systems by maximizing entropy under constraints.

In physics, the missing information on the concrete state of a system is related to the entropy of the system. Entropy is also an elementary information concept. In 1948 Shannon (5), applies Boltzmann's statistical entropy to information theory and gives a new meaning to entropy, relating entropy with the absence/presence of information in a given message. Shannon's entropy is a measure of uncertainty about a system, which can characterize the amount of information. The entropy is defined as our missing information on the actual occurrence of events, given that we only know the probability distribution of the events. As a theoretical foundation, Shannon entropy plays a key role in theory and applications of information theory. Entropy is not only an important concept and a physical quantity in physics and information theory, but also has been widely applied in system theory, cybernetics, economics, management sciences, engineering sciences and other fields.

An axiomatic system is a formal description of a way to establish the scientific truth that flows from a fixed set of assumptions. Shannon and Khinchin proposed an axiomatic system of information entropy, namely, the Shannon-Khinchin (SK) axioms[†] (5). Khinchin had rigorously proven that the Shannon entropy is the only function that satisfies all four of the Khinchin axioms (6). That is to say, the SK axioms are equivalent to Shannon entropy function, and Shannon information function is uniquely determined by all four SK axioms. This axiomatic system ensures the uniqueness of the classical Boltzmann entropy or Shannon entropy form, $H_{BGS} = -k\sum_{i=1}^{W} p_i \log p_i$. However, it is well known that there are some limitations of Shannon entropy, which have restricted its potential applications in physics including complex systems with long-range interactions, scattering processes, and hydrodynamic turbulence (1–4). Daroczy and Tsallis generalized Shannon entropy with adjustable parameters (2, 7), respectively. Therefore, there are other forms of information measures besides Shannon entropy. Beck(8) discussed some of the most important examples of information measures that are useful for the description of complex systems.

To characterize strongly interacting statistical systems within a thermodynamical framework, in particular, complex systems, a series of generalized entropies have been proposed in the past. Based on Shannon-Khinchin axiomatics, Santos (9) and Abe (10) proposed a set of axioms which are appropriate to the Tsallis entropy function in a broader sense, respectively. Santos-Abe axioms is as follows: (i) the entropy is a continuous function of the probabilities $\{p_i\}$; (ii) the entropy is a monotonically increasing function of the number of states $N$ in the case of equidistribution $p_i = 1/N$; (iii) the entropy satisfies the pseudoadditivity relation $H_\beta(A+B)/k = H_\beta(A)/k + H_\beta(B)/k + (1-q)H_\beta(A)H_\beta(B)/k^2$ ($A$ and $B$ being two independent systems, $\beta$ is a constant and $k$ is a positive constant), and (iv) the entropy satisfies the relation $H_\beta(\{p_i\}) = H_\beta(p_L, p_M) + p_L^\beta H_\beta(\{p_i/p_L\}) + p_M^\beta H_\beta(\{p_i/p_M\})$, where $p_L + p_M = 1$ ($p_L = \sum_{i=1}^{N_L} p_i$ and $p_M = \sum_{i=N_L+1}^{N} p_i$). They have proved, along Shannon's line, that the unique function that satisfies all these properties is Tsallis entropy $H_\beta^T(p_1, p_2, \cdots, p_N) = \frac{k}{\beta-1}(1 - \sum_{i=1}^{N} p_i^\beta)$, $\beta \neq 1$. That is to say, their axiomatics is equivalent to Tsallis entropy. However, Santos-Abe axioms can still not unify Shannon entropy, Daroczy entropy, Tsallis entropy, Rényi entropy (11) and other forms of information measures into an axiomatics.

---

[†] Shannon–Khinchin axioms: (K1) Entropy is a continuous function of the probabilities $p_i$ only, i.e., $H$ should not explicitly depend on any other parameters. (K2) Entropy is maximal for the equidistribution $p_i = 1/W$; from this, the concavity of $H$ follows. (K3) Adding a state $W+1$ to a system with $p_{W+1} = 0$ does not change the entropy of the system; from this, $H(0) = 0$ follows. (K4) Entropy of a system composed of two subsystems, A and B, is $H(A+B) = H(A) + H(B/A)$.

Hanel and Thurner (12) proved that under the first three Shannon-Khinchin axioms and the condition assuming that entropy has to be of the form $H = \sum_i g(p_i)$ each statistical system is characterized by a unique pair of scaling exponents (c, d) in the large size limit. These exponents define equivalence classes for all interacting and non-interacting systems. However, if the form $H = \sum_i g(p_i)$ is canceled, then this axiomatic framework is too broad to reflect the nature of some complex systems. The main purpose of the paper is to find a set of reasonable entropy axioms, and a new axiomatic framework of entropy is presented as follows.

**Upper entropy axioms**

Suppose that $\Omega = \{(p_1, p_2, \cdots, p_N) | p_i \geq 0, i = 1, 2, \cdots, N, \sum_{i=1}^{N} p_i = 1\}$, $H$ is a mapping from $\Omega$ into $R$ (all real numbers). $H$ is called an upper entropy function if it satisfies the following axiomatic requirements:

(i) Continuity: $H(p_1, p_2, \cdots, p_n, \cdots)$ is a continuous function, and $H(p_1, p_2, \cdots, p_n, \cdots) = 0$ if and only if there exists a positive integer $k$, such that, $p_k = 1$.

(ii) Symmetry: $H(p_1, p_2, \cdots, p_n, \cdots)$ is symmetric for arbitrary $p_i, i = 1, 2, \cdots$.

(iii) Increasing property: if $q_i \geq 0, i = 1, 2, \cdots, m$, and $p_n = \sum_{i=1}^{m} q_i$, then

$$0 \leq H(p_1, p_2, \cdots, p_{n-1}, q_1, q_2, \cdots q_m, p_{n+1}, \cdots) - H(p_1, p_2, \cdots, p_{n-1}, p_n, p_{n+1}, \cdots)$$
$$\leq p_n H(\frac{q_1}{p_n}, \frac{q_2}{p_n}, \cdots, \frac{q_m}{p_n})$$   [1]

(iv) Maximum: for any finite positive integer $N$,

$$H(p_1, p_2, \cdots, p_N) \leq H(\frac{1}{N}, \frac{1}{N}, \cdots, \frac{1}{N}).$$   [2]

Beck(8) pointed out that the entire formalism of statistical mechanics can be regarded as being based on maximizing the entropy (= missing information) of the system under consideration subject to suitable constraints, and hence naturally the question arises how to measure this missing information. In principle general information measures (that contain the Shannon information as a special case) can be chosen. From the view of information theory, we illuminate upper entropy axioms as follows. Axiom (i) shows that entropy is minimum when there is no uncertainty in a system, while entropy is 0 when each event is certain. Since information entropy depicts the average uncertainty before a message is sent, axiom (ii) shows that entropy is only related to the overall statistical properties of information source. If the statistical properties of two information sources contain the same message and probability distribution, then the entropies of the information sources are of the same. Axiom (iii) shows that, assuming that any one element of an information source is divided into $m$ parts, and the sum of the probabilities of $m$ parts is equal to the probability of the original element of the information source, then the entropy of the new information source may not reduce, and the increasing amount may not exceed the right side of formula [1]. Axiom (iv) means the information measure $H$ takes on an absolute maximum for the uniform distribution $(\frac{1}{N}, \frac{1}{N}, \cdots, \frac{1}{N})$, any other probability distribution has an information content that is less than or equal to that of the uniform distribution.

In this axiomatics, besides Shannon entropy, other functions may also be used to measure information uncertainty.

## Upper entropy function properties

(1) **Expansibility**

For any given $n$, then

$$H(p_1, p_2, \cdots, p_{n-1}, p_n, \cdots) = H(p_1, p_2, \cdots, p_{n-1}, 0, p_n, \cdots).\qquad [3]$$

**Proof**: Axiom (iii) leads to

$$H(p_1, p_2, \cdots, p_{n-1}, p_n, p_{n+1}, \cdots) \leq H(p_1, p_2, \cdots, p_{n-1}, 0, p_n, p_{n+1}, \cdots)$$
$$\leq H(p_1, p_2, \cdots, p_n, \cdots) + p_n H(0,1)$$

From axiom (i), we have $H(0,1) = 0$, thus

$$H(p_1, p_2, \cdots, p_{n-1}, p_n, p_{n+1}, \cdots) \leq H(p_1, p_2, \cdots, p_{n-1}, 0, p_n, p_{n+1}, \cdots) \leq H(p_1, p_2, \cdots, p_n, \cdots),$$

and hence

$$H(p_1, p_2, \cdots, p_{n-1}, p_n, p_{n+1}, \cdots) = H(p_1, p_2, \cdots, p_{n-1}, 0, p_n, p_{n+1}, \cdots)$$

The expansibility means the information measure $H$ should not change if the sample set of events is enlarged by another event that has probability zero.

(2) **Upper subadditivity**

If $X$ and $Y$ are two independent information sources, then

$$H(X+Y) \leq H(X) + H(Y)$$

Namely, let us denote $X = (p_1, p_2, \cdots, p_n)$ and $Y = (q_1, q_2, \cdots q_m)$ as two probability distributions, it follows immediately that

$$H(p_1 q_1, p_1 q_2, \cdots p_1 q_m, p_2 q_1, p_2 q_2, \cdots p_2 q_m, \cdots, p_n q_1, p_n q_2, \cdots p_n q_m)$$
$$\leq H(p_1, p_2, \cdots, p_n) + H(q_1, q_2, \cdots, q_m).\qquad [4]$$

**Proof**: From axiom (iii), we have

$$H(p_1 q_1, p_1 q_2, \cdots p_1 q_m, p_2 q_1, p_2 q_2, \cdots p_2 q_m, \cdots, p_n q_1, p_n q_2, \cdots p_n q_m)$$
$$\leq H(p_1, p_2 q_1, p_2 q_2, \cdots p_2 q_m, \cdots, p_n q_1, p_n q_2, \cdots p_n q_m)$$
$$+ p_1 H(q_1, q_2, \cdots q_m)$$
$$\leq H(p_1, p_2, \cdots, p_n q_1, p_n q_2, \cdots p_n q_m)$$
$$+ p_2 H(q_1, q_2, \cdots q_m) + p_1 H(q_1, q_2, \cdots q_m)$$
$$\leq H(p_1, p_2, \cdots, p_n)$$
$$+ p_n H(q_1, q_2, \cdots q_m) + \cdots + p_2 H(q_1, q_2, \cdots q_m) + p_1 H(q_1, q_2, \cdots q_m)$$
$$= H(p_1, p_2, \cdots, p_n) + H(q_1, q_2, \cdots q_m)$$

(3) **Upper strong subadditivity**

For information sources $X$ and $Y$, then

$$H(X+Y) \leq H(X) + H(Y/X)$$

Namely, let us denote $X = (p_1, p_2, \cdots, p_n)$ and $Y = (q_1, q_2, \cdots q_m)$ as two probability distributions, the conditional probability is used to describe the relationship between them.

$$P(Y = y_j \mid X = x_i) = p_{ij}, \qquad i = 1, 2, \cdots, n, \; j = 1, 2, \cdots m$$

Writing the entropy of $p_{ij}$ as $H(Y/X)$, the conditional entropy is given by

$$H(Y/X) = \sum_{i=1}^{n} p_i H(p_{i1}, p_{i2}, \cdots, p_{im}),\qquad [5]$$

Then

$$H(p_1 p_{11}, p_1 p_{12}, \cdots p_1 p_{1m}, p_2 p_{21}, p_2 p_{22}, \cdots p_2 p_{2m}, \cdots, p_n p_{n1}, p_n p_{n2}, \cdots p_n p_{nm})$$
$$\leq H(p_1, p_2, \cdots, p_n) + \sum_{i=1}^{n} p_i H(p_{i1}, p_{i2}, \cdots, p_{im}).\qquad [6]$$

The proof is similar to that of upper subadditivity (omitted), upper strong subadditivity is easily proved.

(4) **Extremal monotonicity**

For a given set of $N$ equiprobable states, i.e., $p_i = 1/N$, $H$ is a monotonic increasing function of $N$.

**Proof**: From the expansibility, we have

$$H(\frac{1}{N}, \frac{1}{N}, \cdots, \frac{1}{N}) = H(\frac{1}{N}, \frac{1}{N}, \cdots, \frac{1}{N}, 0),\qquad [7]$$

Further from the maximum, we get

$$H(\frac{1}{N},\frac{1}{N},\cdots,\frac{1}{N},0) \leq H(\frac{1}{N+1},\frac{1}{N+1},\cdots,\frac{1}{N+1},\frac{1}{N+1}),  \quad [8]$$

Substituting Eq. 7 into Eq. 8, we can rewrite Eq. 8 by

$$H(\frac{1}{N},\frac{1}{N},\cdots,\frac{1}{N}) \leq H(\frac{1}{N+1},\frac{1}{N+1},\cdots,\frac{1}{N+1},\frac{1}{N+1}).$$

Therefore, for a given set of $N$ equiprobable states, i.e., $p_i = 1/N$, $H$ is a monotonic increasing function of $N$.

Tsallis (1988) defined a type of entropy, namely Tsallis entropy

$$H_\beta^T(p_1, p_2, \cdots, p_n) = \frac{1}{\beta - 1}(1 - \sum_{i=1}^n p_i^\beta), \quad \beta \neq 1 \quad . \quad [9]$$

When $\beta > 1$, we prove the Tsallis entropy satisfies the upper entropy axioms. Apparently, $H_\beta^T(p_1, p_2, \cdots, p_n)$ satisfies axiom (ii) and (iv). Firstly, we prove (i), If there exists $p_i = 1$, obviously, $H_\beta^T(p_1, p_2, \cdots, p_n) = \frac{1}{\beta - 1}(1 - \sum_{i=1}^n p_i^\beta) = 0$. If $H_\beta^T(p_1, p_2, \cdots, p_n) = \frac{1}{\beta - 1}(1 - \sum_{i=1}^n p_i^\beta) = 0$, there exits $i$, such that, $p_i = 1$, $p_j = 0$, $j = 1, 2, \cdots, i-1, i+1, \cdots, n$. Else, $0 < p_j < 1$, $j = 1, 2, \cdots, i-1, i, i+1, \cdots, n$, thus $p_j^\beta < p_j$, $j = 1, 2, \cdots, n$, $1 = \sum_j p_j^\beta < \sum_j p_j = 1$, which leads to a contradiction by itself. Therefore, axiom (i) is proved.

Secondly, we prove the increasing property. If $q_i \geq 0, i = 1, 2, \cdots, m$, and $p_n = \sum_{i=1}^m q_i$, when $\beta > 1$, then $\frac{1}{\beta - 1} > 0$,

$$-\sum_i^m q_i^\beta = -p_n^\beta + p_n^\beta(1 - \sum_{i=1}^m \frac{q_i^\beta}{p_n^\beta}) \leq -p_n^\beta + p_n(1 - \sum_{i=1}^m \frac{q_i^\beta}{p_n^\beta})$$

Eq. 1 is proved. Therefore, axiom (iii) is proved.

Daroczy(1970) defined Daroczy entropy as follows.

$$H_\beta(p_1, p_2, \cdots, p_n) = \frac{1}{1 - 2^{1-\beta}}(1 - \sum_{i=1}^n p_i^\beta), \quad \beta > 0, \beta \neq 1$$

In analogy with Tsallis entropy, we can derive that Daroczy entropy satisfies upper entropy axioms when $\beta > 1$.

Another type of entropy is defined, namely

$$H_\beta(p_1, p_2, \cdots, p_n) = \frac{1}{\beta - 1}(1 - \sum_{i=1}^n (p_i^\beta - (\beta - 1)p_i \ln p_i)), \quad \beta \neq 1 \quad [10]$$

In analogy with Tsallis entropy, we can also prove that this entropy satisfies upper entropy axioms when $\beta > 1$.

## Lower entropy axioms

Suppose $\Omega = \{(p_1, p_2, \cdots, p_N) | p_i \geq 0, i = 1, 2, \cdots, N, \sum_{i=1}^N p_i = 1\}$, $H$ is a mapping from $\Omega$ into $R$ (all real numbers). $H$ is called a lower entropy function if it satisfies the following four entropy axioms.

(i) Continuity: the property is the same as the continuity in upper entropy axioms.

(ii) Symmetric expansibility: $H(p_1, p_2, \cdots, p_n, \cdots)$ is symmetrical for arbitrary $p_i$, and

$$H(p_1, p_2, \cdots, p_{n-1}, p_n, \cdots) = H(p_1, p_2, \cdots, p_{n-1}, 0, p_n, \cdots), \quad [11]$$

(iii) Increasing property: if $q_i \geq 0, i = 1, 2, \cdots, m$, and $p_n = \sum_{i=1}^m q_i$, then

$$H(p_1, p_2, \cdots, p_{n-1}, q_1, q_2, \cdots q_m, p_{n+1}, \cdots) - H(p_1, p_2, \cdots, p_n, \cdots) \geq p_n H(\frac{q_1}{p_n}, \frac{q_2}{p_n}, \cdots, \frac{q_m}{p_n}). \quad [12]$$

(iv) Maximum: the property is the same as the maximum in upper entropy axioms.

When $\beta < 1$, we prove Tsallis entropy satisfies lower entropy axioms. Evidentlly, $H_\beta^T(p_1, p_2, \cdots, p_n)$ satisfies axiom (ii) and (iv). Firstly, we prove (i). If exists $p_i = 1$, obviously, $H_\beta^T(p_1, p_2, \cdots, p_n) = \frac{1}{\beta - 1}(1 - \sum_{i=1}^n p_i^\beta) = 0$. If $H_\beta^T(p_1, p_2, \cdots, p_n) = \frac{1}{\beta - 1}(1 - \sum_{i=1}^n p_i^\beta) = 0$, there exits $i$, such that, $p_i = 1$, $p_j = 0$, $j = 1, 2, \cdots, i-1, i+1, \cdots, n$. Else, $0 < p_j < 1$, $j = 1, 2, \cdots, i-1, i, i+1, \cdots, n$, thus, $p_j^\beta > p_j$, $j = 1, 2, \cdots, n$, $1 = \sum_j p_j^\beta > \sum_j p_j = 1$, which leads to a contradiction by itself. Therefore, nonnegativity axiom is proved.

Secondly, we prove the increasing property. If $q_i \geq 0, i = 1, 2, \cdots, m$, and $p_n = \sum_{i=1}^m q_i$. Since $1 - \sum_{i=1}^m \frac{q_i^\beta}{p_n^\beta} \leq 0$, and $p_n^\beta \geq p_n$, thus,

$$-\sum_i^m q_i^\beta = -p_n^\beta + p_n^\beta(1 - \sum_{i=1}^m \frac{q_i^\beta}{p_n^\beta}) \leq -p_n^\beta + p_n(1 - \sum_{i=1}^m \frac{q_i^\beta}{p_n^\beta})$$

when $\beta < 1$, $\frac{1}{\beta - 1} < 0$, then, Eq. 12 is proved. Therefore, axiom (iii) is proved.

Similar to that of Tsallis entropy, we can prove that Daroczy entropy satisfies lower entropy axioms when $0 < \beta < 1$.

Similar to that of Tsallis entropy, we can also prove that entropy (Eq.10) satisfies lower entropy axioms when $0 < \beta < 1$

**Theorem 1**. If a function satisfies both upper entropy axioms and lower entropy axioms, then it is a Shannon entropy.

**Theorem 2 (lower strong subadditivity).** If $H$ is a lower entropy function, for systems (information sources) $X$ and $Y$, then

$$H(X + Y) \geq H(X) + H(Y/X). \quad [13]$$

Abe (1997) introduced a kind of entropy, namely Abe entropy (8)

$$H_\beta^{Abe}(p_1, p_2, \cdots, p_n) = -\frac{1}{\beta - \beta^{-1}} \sum_{i=1}^n (p_i^\beta - p_i^{-\beta}), \quad \beta \neq 1. \quad [14]$$

From Eq. 14, we have

$$H_\beta^{Abe}(p_1, p_2, \cdots, p_n) = \frac{\beta}{\beta + 1} H_\beta^T(p_1, p_2, \cdots, p_n) + \frac{1}{\beta + 1} H_{\beta^{-1}}^T(p_1, p_2, \cdots, p_n),$$

which implies that Abe entropy is not either an upper entropy function or a lower entropy function.

We may replace axiom (iii) in upper entropy axioms or lower entropy axioms by the less stringent condition, which just states that the entropy of independent systems should be additive. In this case one ends up with other information measures which are called the Rényi entropies (8). These are defined as

$$H_\beta^R(p_1, p_2, \cdots, p_n) = -\frac{1}{\beta - 1} \ln \sum_{i=1}^n p_i^\beta, \quad \beta \neq 1. \quad [15]$$

However, Rényi entropies satisfy neither upper entropy axioms nor lower entropy axioms.
Therefore, these requirements of the upper entropy axioms and the lower entropy axioms are also stronger. In order to find a fairly simple axiomatics, we will develop weak upper entropy axioms and weak lower entropy axioms as follows.

**Weak upper entropy axioms**

Axioms (i), (ii) and axiom (iv) in lower entropy axioms are kept, and axiom (iii) is replaced by the following more general version: new axiom (iii), if system (information source) $X$ and system (information source) $Y$ are independent, then

$$H(X + Y) \leq H(X) + H(Y). \quad [16]$$

**Weak lower entropy axioms**

Axioms (i), (ii) and axiom (iv) in lower entropy axioms are kept, and axiom (iii) is replaced by the following more general version: new axiom (iii), if system (information source) $X$ and system (information source) $Y$ are independent, then

$$H(X+Y) \geq H(X) + H(Y) \:. \qquad [17]$$

Rényi entropy is not only a weak upper entropy function, but also a weak lower entropy function.

Landsberg-Vedral entropy is defined as follows

$$H_\beta^L(p_1, p_2, \cdots, p_n) = -\frac{1}{\beta-1}\left(\frac{1}{\sum_{i=1}^n p_i^\beta} - 1\right), \quad \beta \neq 1 \:. \qquad [18]$$

When $\beta > 1$, Landsberg-Vedral entropy is a weak lower entropy function. When $\beta < 1$, Landsberg-Vedral entropy is a weak upper entropy function.

The following entropy (12)

$$H_\gamma(p_1, p_2, \cdots, p_n) = \sum_i p_i \ln^{1/\gamma}(1/p_i), \quad \gamma > 0 \:. \qquad [19]$$

When $\gamma \geq 1$, the entropy is a weak upper entropy function. When $0 < \gamma \leq 1$, the entropy is a weak lower entropy function.

The following entropy (12)

$$H_\beta(p_1, p_2, \cdots, p_n) = \sum_i p_i^\beta \ln(1/p_i), \quad \beta > 0 \:. \qquad [20]$$

When $0 < \beta \leq 1$, it is a weak lower entropy function. When $\beta \geq 1$, it is a weak upper entropy function.

**Discussion**

Shannon entropy plays a key role in the development and applications of information theory, and also has many applications in physics, management sciences, system sciences and so on. Because it still has many limitations, different forms of information measure such as Daroczy entropy and Tsallis entropy have been proposed. Tsallis statistics is just one example of many possible new statistics. In general, complex nonequilibrium problems may require different types of superstatistics. However, as axioms of metric spaces, whether there is a more reasonable axiomatic framework for superstatistics is an interesting question. To proceed with our exploration, upper entropy axioms and lower entropy axioms are proposed in the paper. This axiomatics provides a framework, which can integrate Shannon entropy, Daroczy entropy, Tsallis entropy and so on. We hope that the work might contribute to the superstatistics theory.

**ACKNOWLEDGMENTS.** Support this work was provided in part by the Shanghai First-class Academic Discipline Project, China (Grant No. S1201YLXK), and supported by the Hujiang Foundation of China (A14006).

**References**

[1] Beck C, Cohen EGD (2003) Superstatistics. *Physica A* 322: 267 -275.
[2] Tsallis C (1988) Possible generalization of Boltzmann-Gibbs statistics. *Journal of Statistical Physics* 52: 479-487.
[3] Asgarani S (2013) A set of new three-parameter entropies in terms of a generalized incomplete Gamma function. *Physica A* 392: 1972-1976.
[4] Hanel R, Thurner S, Gell-Mann M (2011) Generalized entropies and the transformation group of superstatistics. *Proc Natl Acad Sci USA* 108:6390-6395.
[5] Shannon C E (1948) A mathematical theory of communication. *Bell Syst Tech J* 27:379-423 623-656 .
[6] Khinchin A I (1957) Mathematical Foundations of Information Theory (Dover, New York).
[7] Daróczy Z (1970) Generalized information functions. *Information and Control* 16:36-51.
[8] Beck C (2009) Generalised information and entropy measures in physics. *Contemporary Physics*. 50: 495–510.
[9] Dos Santos, R J (1997) Generalization of Shannon's theorem for Tsallis entropy. *J.Math. Phys.* 38: 4104-4109.
[10] Abe S (2000) Axioms and uniqueness theorem for Tsallis entropy. *Phys.Lett.A.* 271:74-80.
[11] Rényi A (1961) On measures of entropy and information. *In Fourth Berkeley Symposium on Mathematical Statistics and Probability* (University of California Press, Berkeley).
[12] Hanel R, Thurner S (2011) A comprehensive classification of complex statistical systems and an axiomatic derivation of their entropy and distribution functions. Europhysics Letters 93: 20006